\documentclass[10pt,journal]{IEEEtran}

\usepackage{amsmath,cite,xspace}
\usepackage{graphicx}
\usepackage{epsfig}
\usepackage{epstopdf}
\usepackage{bm}
\newcommand{\subparagraph}{}
\usepackage[table,xcdraw]{xcolor}
\usepackage{color}
\usepackage{amssymb}

\newcommand\blfootnote[1]{%
  \begingroup
  \renewcommand\thefootnote{}\footnote{#1}%
  \addtocounter{footnote}{-1}%
  \endgroup
}

\allowdisplaybreaks

\usepackage{mathtools}

\newcommand{\bp}{{\bf p}}

\newcommand{\bx}{{\bf x}}

\newcommand{\bv}{{\bf v}}

\graphicspath{{./Figures/}} \DeclareGraphicsExtensions{.eps,.pdf,.jpg,.mps,.png}

\title{Solutions of DC OPF are Never AC Feasible}

 \author{\IEEEauthorblockN{Kyri Baker}\\ University of Colorado Boulder\\Boulder, CO, USA 80309}
\makeatletter
\let\old@ps@headings\ps@headings
\let\old@ps@IEEEtitlepagestyle\ps@IEEEtitlepagestyle
\def\psccfooter#1{%
    \def\ps@headings{%
        \old@ps@headings%
        \def\@oddfoot{\strut\hfill#1\hfill\strut}%
        \def\@evenfoot{\strut\hfill#1\hfill\strut}%
    }%
    \def\ps@IEEEtitlepagestyle{%
        \old@ps@IEEEtitlepagestyle%
        \def\@oddfoot{\strut\hfill#1\hfill\strut}%
        \def\@evenfoot{\strut\hfill#1\hfill\strut}%
    }%
    \ps@headings%
}

\begin{document}
 
\maketitle
\begin{abstract}
In this paper, we analyze the relationship between generation dispatch solutions produced by the DC optimal power flow (DC OPF) problem by the AC optimal power flow (AC OPF) problem. While there has been much previous work in analyzing the approximation error of the DC assumption, the AC feasibility has not been fully explored, although difficulty achieving AC feasibility is known in practice. Here, we consider the set of feasible points in a standard DC OPF problem and the set of feasible points in a standard AC OPF problem. Under some very light assumptions, we show that the intersection of these sets is the empty set; i.e., that no solution to the DC OPF problem will satisfy the AC power flow constraints. Then, it is demonstrated that even with generation adjustments in DC OPF to account for losses, DC OPF solutions are still not AC feasible.
\end{abstract}
\blfootnote{K. Baker is an Assistant Professor in the Department of Civil, Enviornmental, and Architectural Engineering at the University of Colorado Boulder, Boulder, CO, USA. Email: kyri.baker@colorado.edu}
\vspace*{-5mm}
\section{Introduction}

The DC optimal power flow (DC OPF) approximation of the AC optimal power flow (AC OPF) problem is widely used in current power system operation due to its convexity and computational benefits that allow it to be solved on fast timescales. However, it is well-known that the DC OPF can provide, in some cases, a poor approximation of actual AC power flows and  locational marginal prices (LMPs), resulting in physically unrealizable system states. A report from FERC states that many grid operators calculate DC OPF setpoints and modify these until they are AC feasible in a ``DC OPF with AC Feasibility" iterative procedure \cite{FERC11}.

There has been much previous work in analyzing the approximation error between DC OPF and AC OPF solutions \cite{DanErrorBounds16, OPFApprox17, DCrevisit09, CoffrinApprox12}. For example, DC OPF is often used to calculate market prices, and previous works have looked at the difference in LMPs \cite{Overbye04} resulting from DC OPF and AC OPF. Outside of analyzing the differences in objective function value and approximation error between DC and AC OPF, it is important to consider the physical feasibility of the DC OPF solution. Namely, under which circumstances does the solution to a DC OPF problem produce a physically realizable solution (i.e., one that satisfies the AC power flow equations)? The difficulty of obtaining an AC feasible solution from DC OPF is well known; for example, in \cite{lowinfeasible} the authors state ``a solution of DC OPF may not be feasible," and in \cite{dmitry19} the authors state that the DC OPF solution is ``typically AC infeasible." In this note, we show that these statements can actually be made much stronger, because the solution to the standard DC OPF problem is actually always AC infeasible. Considering the origin of these differences may help us gain insight into how to address the AC feasibility of DC OPF.

We show that under some light assumptions (at least one power flow is present in the network, and loads are modeled as constant (P,Q) loads, both of which are typical assumptions in OPF), the set of points within the feasible region of the DC OPF problem and the set of points within the feasible region of the AC OPF problem have an empty intersection. While this has been observed in practice, we demonstrate mathematically why the DC OPF is never AC feasible. In addition, we show that even if a slack bus is used to make up for the lack of losses in DC OPF, the DC OPF voltages and angles are still not AC feasible. This has implications for future power system operation, considering many system operators still use DC OPF and are thus always required to go through an additional iterative procedure to ensure AC feasibility \cite{FERC11}. As shown here, as system loading increases, the gap between the DC and AC OPF solutions also increases. 


Towards the aforementioned issues, this paper offers two main contributions:

\begin{itemize}
    \item We present what we believe is the first formal analysis demonstrating there is no overlap between the DC OPF and AC OPF feasible regions. The benefit of performing a mathematical analysis here rather than further relying on heuristics and physical observations is that it offers us insights into origins of the feasibility gap.
    \item We further show that even if the total generation in AC and DC OPF is made equal by inclusion of fictitious nodal demand or another technique, the solution obtained by DC OPF does not satisfy the AC power flow equations.
\end{itemize}

\section{Optimal Dispatch and Power Flow}
To provide the necessary background to discuss the feasibility of the DC OPF solutions in the AC OPF problem, we first briefly summarize these problems which seek to optimize generation dispatch in transmission networks. First, define coefficients $a_j$, $b_j$, and $c_j$ as the operational costs associated with generator $j$. Let set $\cal{G}$ be the set of all generators in the network, $\cal{N}$ be the set of all nodes (buses) in the network, $\cal{L}$ be the set of all lines (branches) in the network, and $\mathcal{G}_i$ be the set of generators connected to bus $i$. Define $p_{l,j}$ ($q_{l,j}$) as the total active (reactive) power consumption at node $j$, $p_{g,j}$ ($q_{g,j}$) is the active (reactive) power output of generator $j$, and $\underline{p}_{g,j}$ ($\underline{q}_{g,j}$) and $\overline{p}_{g,j}$ ($\overline{q}_{g,j}$) are lower and upper limits on active (reactive) power generation, respectively. The complex voltage at node $i$ has magnitude $|v_i|$ and phase angle $\theta_i$, and the difference in phase angle between neighboring buses $i$ and $m$ is written as $\theta_{im}$. Values $G_{im}$ and $B_{im}$ are real and reactive parts of entry $(i,m)$ in the admittance matrix, respectively.




\subsection{AC Optimal Power Flow}

The AC Optimal Power Flow (AC OPF) model is typically considered the ground truth for estimating physical power flows throughout the network, as it includes line losses, network parameters, and reactive power. It can be written with the AC power flow equations in polar form as follows:

\begin{subequations} \label{eqn:ACOPF}
\begin{align} 
\min_{\substack{\bv, \bp_g}} &\sum_{j \in \cal{G}} a_j p_{g,j}^2 + b_j p_{g,j} + c_j \\
\mathrm{s.t:~}\nonumber\\
&\hspace{-0.8cm}|v_i| \sum_{m \in \mathcal{N}} |v_m|(G_{im}\cos(\theta_{im}) + B_{im}\sin(\theta_{im})) \nonumber\\
&\hspace{2cm}=~ p_{l,i} - \sum_{k \in \mathcal{G}_i} p_{g,k}, ~~\forall i \in \cal{N}\label{eqn:PF_constr}\\
&\hspace{-0.8cm}|v_i| \sum_{m \in \mathcal{N}} |v_m|(G_{im}\sin(\theta_{im}
) - B_{im}\cos(\theta_{im}))\nonumber \\
&\hspace{2cm}=~ q_{l,i} - \sum_{k \in \mathcal{G}_i} q_{g,k}, ~~\forall i \in \cal{N}\label{eqn:PF_constr2}\\
&\underline{p}_{g,j} \leq p_{g,j} \leq \overline{p}_{g,j}, ~~\forall j \in \cal{G} \\
&\underline{q}_{g,j} \leq q_{g,j} \leq \overline{q}_{g,j}, ~~\forall j \in \cal{G} \label{eqn:Q_constr}\\
&\underline{|v|} \leq |v_i| \leq \overline{|v|}, ~~\forall i \in \cal{N}. \label{eqn:V_constr}
 \end{align}
\end{subequations}

\noindent where $\bp_g$ is a vector comprising the active power generation $p_{g,j}$ at each generator $j \in \cal{G}$ and $\bv$ is $2n$-dimensional vector comprising the unknown voltage magnitudes and angles. It is clear from \eqref{eqn:PF_constr} and \eqref{eqn:PF_constr2} that the overall optimization problem is nonconvex. Note that here we omit line flow constraints, but these can be modeled as constraints on apparent power flows, line currents, or voltage angle differences \cite{lineFlows13}.

\subsection{DC Optimal Power Flow}
The constraints within the DC Optimal Power Flow problem (DC OPF) are linear approximations of the actual nonlinear AC power flows. The DC approximation is derived from multiple physical assumptions and observations. First, in transmission networks, the line resistance $R_{im}$ is typically significantly less than the line reactance $X_{im}$; thus, $B_{im}$ can be approximated to $-\frac{1}{X_{im}}$. Second, the phase angle difference between any two buses is typically small and usually does not exceed $30^{\circ}$. From this, we can use the small angle approximation to approximate $\sin(\theta_{im}) \approx \theta_{im}$ and $\cos(\theta_{im}) \approx 1$. Third, transmission level voltage magnitudes are typically very close to $1.0$ p.u. during normal operation. 

Lastly, from these initial assumptions, and confirmed by the fact that reactive power is a localized phenomenon that cannot travel long distances, we see that the magnitude of the reactive power flow on the lines (denote this as $Q_{im}$ for line $im$) is significantly less than the magnitude of the active power flow on the lines (denote this as $P_{im}$ for line $im$). Using these assumptions when studying the AC power flow equations, equations \eqref{eqn:PF_constr} and \eqref{eqn:PF_constr2} simplify and we are left with the following DC OPF problem. Note that in some formulations, the objective can also be linear instead of quadratic.

\begin{subequations} \label{eqn:DCOPF}
\begin{align} 
\min_{\substack{\bp_g}} &\sum_{j \in \cal{G}} a_j p_{g,j}^2 + b_j p_{g,j} + c_j \\
\mathrm{s.t:~}&p_{l,i} - \sum_{k \in \mathcal{G}_i} p_{g,k} = \sum_{m \in \mathcal{N}} B_{im}\theta_{im}, ~~\forall i \in \cal{N}\label{eqn:DCPowerBalance}\\
&\hspace{.5cm}-F_{im} \leq B_{im}\theta_{im} \leq F_{im}, ~~\forall im \in \cal{L}\label{eqn:DClineflow}
\end{align}
\end{subequations}

\noindent where $F_{im}$ represents the limit on the magnitude of the line flows on line $im$. Note that physically, transmission line flows are limited by the amount of current that can safely flow through the line. We can write the current flow limit on line $im$, $|I_{im}|$, in terms of the complex power $S_{im}$, real power $P_{im}$, and reactive power $Q_{im}$ flowing from bus $i$ to bus $m$ and the complex voltage at bus $i$:
\begin{align*}
|I_{im}|= \Big|\Big(\frac{S_{im}}{v_i}\Big)^*\Big| = \Big(\frac{\sqrt{(P_{im}^2+Q_{im}^2)}}{|v_i|}\Big),
\end{align*}

\noindent and by using the DC approximations stated above, namely that $P_{im} >> Q_{im}$ and $|v_i| = 1.0$ p.u., we can write the current flow limit in terms of power

\begin{align*}
|I_{im}| \approx \sqrt{P_{im}^2} \approx F_{im},
\end{align*}

\noindent which justifies \eqref{eqn:DClineflow}. More information on the derivation of the DC approximation can be found in \cite{LineLossDC}.

\section{AC Feasibility}
Towards determining the AC feasibility of DC OPF, define the set of feasible $\bp_g$ satisfying \eqref{eqn:DCOPF} as $\cal{Y_{DC}} \in \mathbb{R}^{|\cal{G}|}$, and the set of feasible $\bp_g$ satisfying \eqref{eqn:ACOPF} as $\cal{Y_{AC}} \in \mathbb{R}^{|\cal{G}|}$. Note that for the first part of the feasibility analysis (without loss adjustments), it is sufficient to consider generation only. Next, we state the assumptions used in the analysis and results in this paper, which are generally reasonable assumptions made in most AC OPF analyses. \vspace{1mm}

\noindent \emph{Assumption 1:} Assume that the loads in the network are modeled as constant (P, Q) loads \cite{PFNewton}. \vspace{1mm}\\
\noindent \emph{Assumption 2:} Assume that power is flowing on at least one line in the network; i.e., for all feasible solutions $\bx \in \cal{Y_{AC}}$, $\sum_{j \in \cal{G}} p_{g,j} > \sum_{i=1}^{\cal{N}} p_{l,i}$ due to line losses. \vspace{1mm}\\
\noindent \emph{Assumption 3:} Assume the admittance matrix of the system is symmetric.\vspace{1mm} \\
\noindent \emph{Assumption 4:} Assume line resistances and reactances are positive.

From these, we show that $\cal{Y_{DC}} \cap \cal{Y_{AC}} = \emptyset$. Note that because DC OPF neglects reactive power flows, it is not as meaningful to analyze AC feasibility by using reactive power balance equations; instead, active power balance equations are used.

\subsection{AC Feasibility of Economic Dispatch}
It is first intuitive to discover that the solution to the economic dispatch (ED) problem is never AC feasible in practice due to the lack of network flows and line losses. Call $\cal{Y_{ED}}$ the set of feasible solutions to the ED problem. Then define $p_i := |v_i| \sum_{m \in \mathcal{N}} |v_m|(G_{im}\cos(\theta_{im}) + B_{im}\sin(\theta_{im}))$; e.g., define $p_i$ as the left hand side of \eqref{eqn:PF_constr}. If we sum over all loads and generators in (2b), we obtain $\sum_{i \in \cal{N}} p_i = \sum_{i \in \cal{N}} p_{l,i} - \sum_{j \in \cal{G}} p_{g,j}$. From ED, we know that $\sum_{j \in \cal{G}} p_{g,j} = \sum_{i \in \cal{N}} p_{l,i}$. Thus, $\sum_{i \in \cal{N}} p_i = 0$. However, Assumption 2 states that due to line losses, $\sum_{j \in \cal{G}} p_{g,j} > \sum_{i \in \cal{N}} p_{l,i}$. If this must be true, we obtain $\sum_{i \in \cal{N}} p_i \neq 0$, meaning we have a contradiction and there exists no solution for $\bp_g$ that is an element of both $\cal{Y_{ED}}$ and $\cal{Y_{AC}}$; thus $\cal{Y_{ED}} \cap \cal{Y_{AC}} = \emptyset$.


\subsection{AC Feasibility of DC OPF}
Next, consider the DC OPF problem as outlined in \eqref{eqn:DCOPF}. The difference between ED and DC OPF is the consideration of line flows. However, notice that due to the lack of line losses in DC OPF, there exists a similarity between these two problems. In particular, note that the ED constraint $\sum_{j \in \cal{G}} p_{g,j} = \sum_{i \in \cal{N}} p_{l,i}$ also holds for DC OPF. To illustrate this, consider constraint \eqref{eqn:DCPowerBalance}. Summing the left-hand side and right-hand side over all $i \in \cal{N}$ yields

\begin{align}
\sum_{i \in \cal{N}}p_{l,i} - \sum_{j \in \cal{G}}p_{g,j} = \sum_{i \in \cal{N}} \sum_{m \in \mathcal{N}} B_{im}\theta_{im}. \label{eqn:DCsum}
\end{align}

\noindent Since $B_{im} = B_{mi}$ and $\theta_{im} = -\theta_{mi}$, we are left with

\begin{align}
\sum_{i \in \cal{N}}p_{l,i} - \sum_{j \in \cal{G}}p_{g,j} = 0. \label{eqn:DCsum2}
\end{align}

Because of Assumption 2 and the presence of line losses,  $\sum_{j \in \cal{G}} p_{g,j}$ will always be strictly greater than $\sum_{i \in \cal{N}} p_{l,i}$ in AC OPF. Losses have previously been identified as one of the main sources of DC approximation errors \cite{DCrevisit09}; they are also the main cause of infeasibility, as shown here. Thus, $\cal{Y_{DC}} \cap \cal{Y_{AC}} = \emptyset$. Note that the generation values alone are enough to show this result, and we make no statements about the differences between voltage phase angles in the AC and DC cases.

\begin{figure}[t!]
    \includegraphics[width=0.5\textwidth]{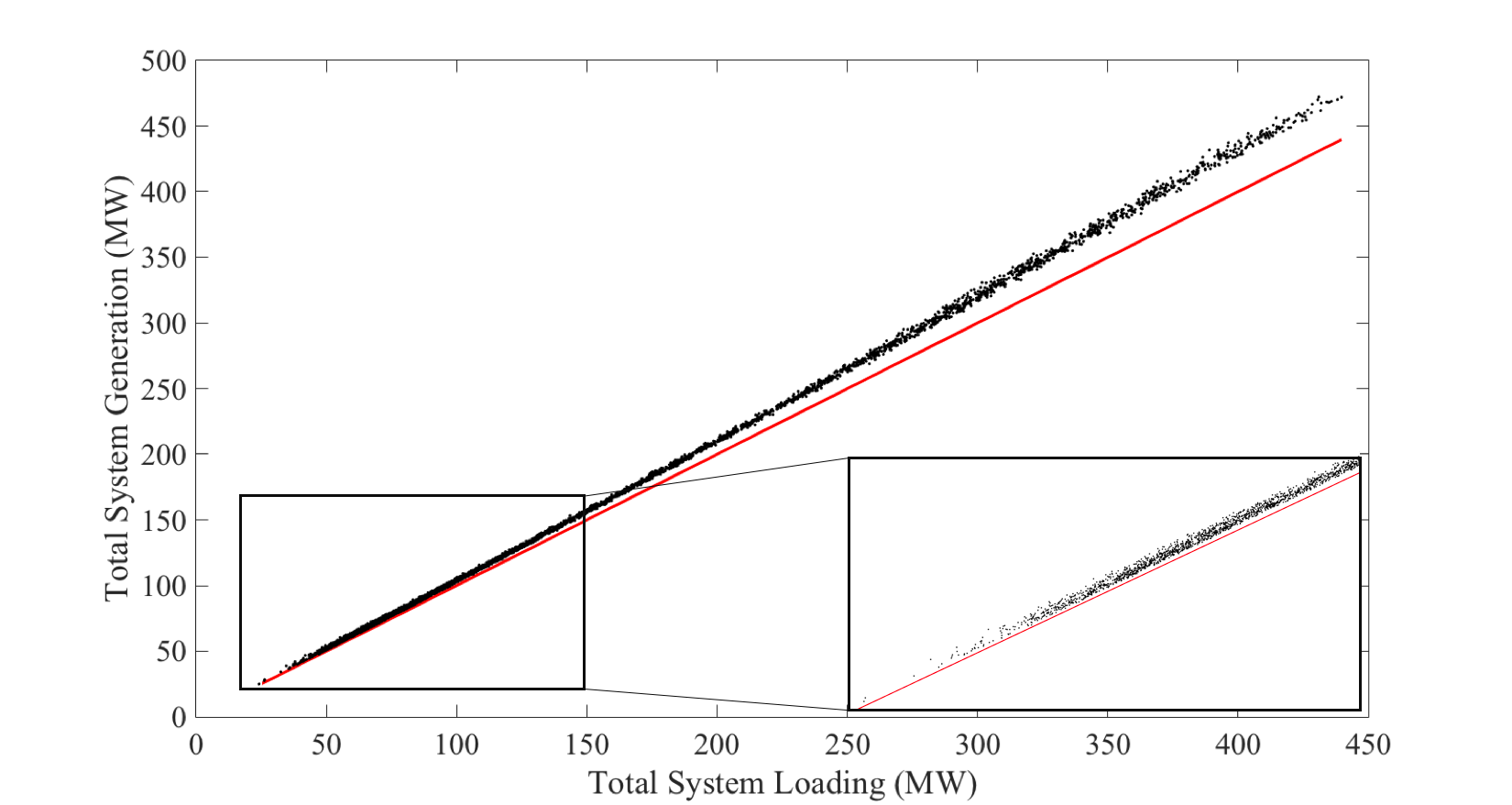}
    \caption{Total generation active power output from DC OPF (red line) and AC OPF (black dots) for varying levels of total system loading. Due to tranmission losses, AC OPF always produces a solution where  $\sum_{j \in \cal{G}} p_{g,j} > \sum_{i \in \cal{N}} p_{l,i}$ whereas DC OPF always produces a solution where  $\sum_{j \in \cal{G}} p_{g,j} = \sum_{i \in \cal{N}} p_{l,i}$. This difference is one of the causes for DC OPF solutions never being AC feasible.}
    \label{fig:totgen}
\end{figure}

To provide a simple illustration of the overall generation in each case, Figure \ref{fig:totgen} shows 5000 feasible runs of both DC and AC OPF for varying levels of random system loading in the IEEE 14-bus system using MATPOWER \cite{MATPOWER}. The red line in the figure shows the total system generation resulting from DC OPF, and the black dots show the total generation resulting from AC OPF. Multiple AC OPF solutions can correspond to one value of total system loading due to the inclusion of losses and distribution of the loading. Notice that this gap increases as system loading and losses increase.


Physically, no dispatch solution is ``AC infeasible." Power flows will always adhere to physics, and if not enough generation is dispatched to meet the load, brownouts or power outages may be experienced. However, the DC OPF solutions are often used in real systems to calculate marginal prices and often in academia to demonstrate algorithmic capabilities, limitations, or to perform research studies. This paper simply calls into question the validity of using DC OPF for these purposes if, in addition to possible optimality differences, feasibility differences also exist.

\section{Feasibility of DC OPF with loss adjustments}
Although often not implemented in research leveraging DC OPF formulations, in reality, of course, losses are not ignored. Generators except the designated slack bus will receive and implement the setpoints determined by DC OPF, and the slack bus will account for the mismatch between the DC OPF-determined generation setpoints and the actual network demand, for example. This will result in the solutions in Fig. 1 overlapping - total generation in the DC OPF will then equal total generation in the AC OPF. Losses can also be included by adding ``fictitious nodal demand" by modifying the load at each bus in order to result in a solution closer to AC OPF \cite{DCloss2, DCloss4}. Now, assuming there is a situation where the generation dispatch solution provided by DC OPF is equal to that provided by AC OPF, we can analyze if the DC OPF solution will be AC feasible, assuming the structure of the DC power flow equations themselves does not change.

\subsection{With small angle approximation}
A DC OPF solution comprised of voltage angles and active power dispatch levels is considered ``AC feasible" if it satisfies the AC power flow equations. We will show with a proof by contradiction that even with the slack bus accounting for system losses or artificially increasing demand, DC OPF results in an AC infeasible point. To that end, consider the active power balance AC power flow equation at some bus $i$ evaluated at a DC OPF solution, where DC OPF necessitates $|v_i| = 1.0$ for all voltage magnitudes in the network. 

\begin{align}
\sum_{m \in \mathcal{N}} (G_{im}\cos(\theta_{im}) + B_{im}\sin(\theta_{im})) =~ p_{l,i} - \sum_{k \in \mathcal{G}_i} p_{g,k}.
\end{align}

Note that the right hand side, by virtue of the DC power flow equations being satisfied, can be rewritten as:

\begin{align}
\sum_{m \in \mathcal{N}} (G_{im}\cos(\theta_{im}) + B_{im}\sin(\theta_{im})) = \sum_{m \in \mathcal{N}} B_{im}\theta_{im}.\label{eq:ACDC1}
\end{align}

\noindent This equality should hold if the DC solution (namely, $|V_i|=|V_m|=1.0$ and $\theta_{im}$, $\forall m \in \mathcal{N}$) are AC feasible. In this case, consider that the small angle approximation holds, and thus

\begin{align}
\sum_{m \in \mathcal{N}} (G_{im} + B_{im}\theta_{im}) = \sum_{m \in \mathcal{N}} B_{im}\theta_{im}.
\end{align}

\noindent Since $G_{im} > 0$, this cannot hold and the solution from DC OPF is not AC feasible.

\subsection{Without small angle approximation}
We can also show the AC infeasibilty without assuming $\theta_{im}$ is small. Two separate cases must be considered under a wider range of angle differences (typically, $-\frac{\pi}{6} < \theta_{im} < \frac{\pi}{6}$ \cite{DCuse05}; here, the range $-\frac{\pi}{2} < \theta_{im} < \frac{\pi}{2}$ is considered).\vspace{2mm}

\noindent \textbf{Case 1:} $0 \leq \theta_{im} < \frac{\pi}{2}$.\\ Rearrange \eqref{eq:ACDC1}:

\begin{align}
\sum_{m \in \mathcal{N}} (G_{im}\cos(\theta_{im})) = \sum_{m \in \mathcal{N}} B_{im}(\theta_{im}-\sin(\theta_{im})).\label{eq:ACDC2}
\end{align}

\noindent Note that within the given range of $\theta_{im}$, the left-hand side must always be strictly positive. In addition, $B_{im} < 0$ and $\theta_{im} > sin(\theta_{im})$ in this region, so the right hand side must be strictly negative. Thus we have a contradiction.\vspace{2mm}

\noindent \textbf{Case 2:} $-\frac{\pi}{2} \leq \theta_{im} < 0$. \\
Differentiate both sides of \eqref{eq:ACDC2} with respect to $\theta_{im}$:

\begin{align}
-\sum_{m \in \mathcal{N}} (G_{im}\sin(\theta_{im})) = \sum_{m \in \mathcal{N}} B_{im}(1-\cos(\theta_{im})). \label{eq:ACDC3}
\end{align}

\noindent Note that in this region, $\sin(\theta_{im}) < 0$, and since $G_{im} > 0$, the left hand side of \eqref{eq:ACDC3} is strictly positive. Since $B_{im} < 0$ and in the given region, $0 < \cos(\theta_{im}) < 1$, the right hand side is strictly negative. Thus, we have a contradiction.

The above analyses can also be applied to the standard DC OPF formulation without including the slack bus accounting for line losses; however, the above is less intuitive than the results in Section III-B.


\section{Improving DC OPF and Conclusions}

In this paper, we mathematically demonstrated how a solution to the DC OPF, under some light assumptions, will never be able to satisfy the AC power flow equations. While this has been observed in practice, and the employment of iterative ``AC feasibility" techniques are currently being used by system operators to modify the original DC OPF problem, it had not yet been shown that a solution satisfying both DC and AC OPF does not exist. Much previous work has been focused on quantifying and measuring the DC approximation error rather than analyzing AC feasibility. The findings in this paper and others indicate that the standard DC OPF may not be a representative problem for grid operators to use, and other linear approximations, or pursuing more efficient ways of solving the AC OPF problem might be more appropriate. Improving line loss modeling has been of particular interest  \cite{DCloss1,DCloss2,DCloss3,CoffrinApprox12}; however, as seen in this paper, much of the infeasibility stems from the structure of the DC power flow equations themselves and the lack of consideration of voltage magnitudes.

It is worth noting that large scale ($>$50,000 bus) AC OPF problems can be solved in minutes using off-the-shelf solvers \cite{Kardos18}. With increases in computing power and availability, AC OPF becomes a more and more tractable problem. The need for post-processing, loss estimation techniques, lack of AC feasibility, and sacrifices in optimality make it challenging to continue to argue for the use of DC OPF. It may be another story for complicated problems such as stochastic OPF, security-constrained OPF, or other mixed-integer or stochastic problems, however. Increasing renewable integration and uncertainty may also expedite needs for accurate models.



\bibliographystyle{IEEEtran}
\bibliography{references.bib}

\end{document}